\documentclass[12pt,twoside]{amsart}
\usepackage{latexsym,amssymb,amsmath}
\numberwithin{equation}{section}
\begin{document}

\title[Symbolic Rees algebras]{Symbolic Rees algebras, vertex covers and 
irreducible representations of Rees cones}

\author{L. A. Dupont}
\address{Departamento de Matem\'aticas\\ 
Centro de Investigaci\'on y de Estudios
\\ Avanzados del IPN\\
Apartado Postal 14--740 \\
07000 Mexico City,
D.F.}
\email{ldupont@math.cinvestav.mx}

\author{R. H. Villarreal}
\address{Departamento de Matem\'aticas\\ 
Centro de Investigaci\'on y de Estudios
\\ Avanzados del IPN\\
Apartado Postal 14--740 \\
07000 Mexico City,
D.F.}
\email{vila@math.cinvestav.mx}

\thanks{Partially supported
by CONACyT grant 49251-F and SNI, M\'exico.}

%\date{1.02.2009}

%Change theorem environments according to your needs...
\theoremstyle{plain}
\newtheorem{theorem}{Theorem}[section]
\newtheorem{lemma}[theorem]{Lemma}
\newtheorem{proposition}[theorem]{Proposition}
\newtheorem{corollary}[theorem]{Corollary}
\theoremstyle{definition}
\newtheorem{definition}[theorem]{Definition} 
\newtheorem{remark}[theorem]{Remark}
\newtheorem{example}[theorem]{Example}



\begin{abstract} Let $G$ be a simple graph and let $I_c(G)$ be its
ideal of vertex 
covers. We give a graph theoretical description of the 
irreducible $b$-vertex covers of $G$, i.e., we describe the minimal 
generators of the symbolic Rees algebra of $I_c(G)$. 
Then we study the irreducible $b$-vertex covers of the blocker of 
$G$, i.e., we study the minimal generators of the 
symbolic Rees algebra of the edge ideal of $G$. We give a
graph theoretical description of the irreducible binary $b$-vertex
covers of the 
blocker of $G$. It is shown that they correspond to irreducible induced 
subgraphs of $G$. As a byproduct we obtain a method, using Hilbert
bases, to obtain all irreducible induced subgraphs of $G$. In
particular we obtain all  
odd holes and antiholes. We study irreducible graphs and give a method 
to construct irreducible $b$-vertex
covers of the blocker of $G$ with high degree 
relative to the number of vertices of $G$.  
\end{abstract}

\maketitle

\section{Introduction}

A {\it clutter\/} $\mathcal C$ with 
vertex set $X=\{x_1,\ldots,x_n\}$ is a family of subsets of $X$,
called edges, none of which is included in another. The set of
vertices and edges of $\mathcal C$ are denoted by $V(\mathcal{C})$ and
$E(\mathcal{C})$ respectively. A basic example  
of a clutter is a graph. Let $R=K[x_1,\ldots,x_n]$ be a polynomial ring 
over a field $K$. The {\it edge ideal\/} of $\mathcal{C}$, 
denoted by $I(\mathcal{C})$, is the ideal of $R$
generated by all monomials $\prod_{x_i\in e}x_i$ such 
that $e\in E(\mathcal{C})$. The assignment $\mathcal{C}\mapsto
I(\mathcal{C})$ 
establishes a natural one to one
correspondence between the family of clutters and the family of 
square-free monomial ideals. Let $\mathcal{C}$ be a clutter and let 
$F=\{x^{v_1},\ldots,x^{v_q}\}$ be the minimal set of generators of its
edge ideal $I=I(\mathcal{C})$. As usual we
use  $x^a$ as an abbreviation for $x_1^{a_1} \cdots x_n^{a_n}$, 
where $a=(a_1,\ldots,a_n)\in \mathbb{N}^n$. The $n\times q$ matrix
with column vectors 
$v_1,\ldots,v_q$ will be denoted by $A$, it is called the {\it incidence
matrix\/} of $\mathcal C$.

The {\it blowup algebra\/} studied here is the {\it symbolic Rees
algebra\/}:
$$
R_s(I)=R\oplus I^{(1)}t\oplus\cdots\oplus I^{(i)}t^i\oplus\cdots\subset R[t],
$$
where $t$ is a new variable and $I^{(i)}$ is the 
$i${\it th\/} symbolic power of $I$. Closely related to $R_s(I)$ is
the {\it Rees algebra\/} of $I$:
$$
R[It]:=R\oplus It\oplus\cdots\oplus I^{i}t^i\oplus\cdots
\subset R[t].
$$

The study of symbolic powers of edge ideals was initiated in 
\cite{ITG} and further elaborated on in
\cite{bahiano,normali,reesclu,clutters,
cover-algebras,sullivant,perfect}. By a result of
Lyubeznik \cite{Lyu3}, $R_s(I)$ is a $K$-algebra of finite type. 
In general the minimal set of generators of $R_s(I)$ as a $K$-algebra
is very hard to describe in terms of $\mathcal C$ (see \cite{bahiano}). 
There are two exceptional cases. If the clutter $\mathcal{C}$ has the
max-flow  
min-cut property, then by a result of \cite{clutters} 
we have $I^i=I^{(i)}$ for all $i\geq 1$, i.e., $R_s(I)=R[It]$. If
${G}$ is a perfect graph, then the minimal generators of
$R_s(I(G))$ are in one to one correspondence with the cliques
(complete subgraphs) of ${G}$ \cite{perfect}. We
shall be interested in 
studying the minimal set of generators of $R_s(I)$ using polyhedral
geometry. Let $G$ be a graph and let $I_c(G)$ be the
Alexander dual of $I(G)$, see definition below. 
Some of the main results of this paper are 
graph theoretical descriptions of the minimal generators of
$R_s(I(G))$ and $R_s(I_c(G))$. In Sections~\ref{blowupcovers1} and
\ref{blowupcovers} we show that both algebras encode combinatorial
information of the graph which can be decoded using integral 
Hilbert bases. 

The {\it Rees cone\/} of $I$, denoted by $\mathbb{R}_+(I)$, is the
polyhedral  
cone consisting  
of the non-negative linear combinations of the set 
$${\mathcal
A}'=\{e_1,\ldots,e_n,(v_1,1),\ldots,(v_q,1)\}\subset\mathbb{R}^{n+1},$$
where $e_i$ is the $i${\it th} unit vector.

A subset $C\subset X$ is called a {\it vertex cover\/} of the clutter
$\mathcal C$ if every edge of $\mathcal C$ contains at least one
vertex of $C$. A subset $C\subset X$ is called a {\it minimal 
vertex cover\/} of the clutter
$\mathcal C$ if $C$ is a vertex cover of $\mathcal{C}$ and no 
proper subset of $C$  is a vertex cover of $\mathcal{C}$.  
Let $\mathfrak{p}_1,\ldots,\mathfrak{p}_s$ be the minimal primes 
of the edge ideal $I=I({\mathcal C})$ and let 
$$
C_k=\{x_i\vert\, x_i\in\mathfrak{p}_k\}\ \ \ \ (k=1,\ldots,s)
$$ 
be the corresponding minimal vertex covers of $\mathcal C$, see 
\cite[Proposition 6.1.16]{monalg}. Recall that the primary 
decomposition of the edge ideal of $\mathcal{C}$ is given by
$$
I(\mathcal{C})=(C_1)\cap (C_2)\cap\cdots\cap (C_s),
$$
where $(C_k)$ denotes the ideal of $R$ generated by $C_k$. 
In particular observe that the height of 
$I(\mathcal{C})$ equals the number of vertices in a minimum vertex
cover of $\mathcal{C}$. This number is called the {\it vertex
covering number\/} of $\mathcal{C}$ and is denoted by
$\alpha_0(\mathcal{C})$. The $i${\it th} {\it symbolic power\/} of
$I$ is given 
by 
\[
I^{(i)}=S^{-1}I^i\cap R\ \mbox{ for }i\geq 1,
\]
where $S=R\setminus\cup_{k=1}^s{\mathfrak p}_i$ and $S^{-1}I^i$ is the
localization of $I^i$ at $S$. In our situation the $i${\it th}
symbolic power of $I$ has a simple expression: 
$$
I^{(i)}=\mathfrak{p}_1^i\cap\cdots\cap \mathfrak{p}_s^i,
$$
see
\cite{monalg}. The Rees cone of $I$ is a finitely generated 
rational cone of dimension $n+1$. Hence by the finite basis theorem
\cite[Theorem~4.11]{webster} there is a 
unique irreducible representation
\begin{equation}\label{okayama-car1} 
{\mathbb R}_+(I)=H_{e_1}^+\cap H_{e_2}^+\cap\cdots\cap
H_{e_{n+1}}^+\cap H_{\ell_1}^+\cap H_{\ell_2}^+\cap\cdots\cap
H_{\ell_r}^+
\end{equation}
such that each $\ell_k$ is in $\mathbb{Z}^{n+1}$, the non-zero
entries of 
each $\ell_k$ are relatively prime, and none of the closed 
halfspaces $H_{e_1}^+,\ldots,
H_{e_{n+1}}^+,H_{\ell_1}^+,\ldots,H_{\ell_r}^+$ can be
omitted from 
the intersection. Here $H_{a}^+$ denotes 
the closed halfspace 
$H_a^+=\{x\vert\, \langle
x,a\rangle\geq 0\}
$
and $H_a$ stands for the hyperplane through the origin with normal
vector $a$, where $\langle\ ,\, \rangle$ denotes the standard 
inner product. The {\it facets\/} (i.e., the proper faces of maximum 
dimension or equivalently the faces of dimension $n$)
of the Rees cone are exactly:
\begin{equation*}
F_i=H_{e_i}\cap {\mathbb R}_+(I), i=1,\ldots,n+1, H_{\ell_1}\cap
{\mathbb R}_+(I), \ldots, 
H_{\ell_r}\cap {\mathbb R}_+(I).
\end{equation*}
According to \cite[Lemma~3.1]{normali} we may always assume 
that $\ell_k=-e_{n+1}+\textstyle\sum_{x_i\in C_k}e_i$ for $1\leq
k\leq s$, i.e., each minimal vertex cover of $\mathcal{C}$ determines
a facet of 
the Rees cone and every facet of the Rees cone satisfying 
$\langle \ell_k,e_{n+1}\rangle=-1$ must be of the
form $\ell_k=-e_{n+1}+\sum_{x_i\in C_k}e_i$ for some minimal vertex
cover $C_k$ of $\mathcal{C}$. This is quite interesting because this
is saying that the Rees cone of $I(\mathcal{C})$ is a carrier of 
combinatorial information of the clutter $\mathcal{C}$. Thus we can
extract the primary decomposition of $I(\mathcal{C})$ from the
irreducible representation of $\mathbb{R}_+(I(\mathcal{C}))$.

Rees cones have been used to study algebraic and combinatorial
properties of blowup algebras of square-free monomial ideals and
clutters  
\cite{normali,reesclu,matrof}. Blowup algebras are interesting 
objects of study in algebra and geometry \cite{Vas}.  

The ideal of {\it vertex covers\/} 
of ${\mathcal C}$ is the square-free monomial ideal
$$
I_c({\mathcal C})=(x^{u_1},\ldots,x^{u_s})\subset R,
$$
where $x^{u_k}=\prod_{x_i\in C_k}x_i$. Often the ideal
$I_c(\mathcal{C})$ is
called the {\it Alexander dual\/} of $I(\mathcal{C})$. The
clutter $\Upsilon({\mathcal C)}$ associated to $I_c({\mathcal C})$ is
called  
the {\it blocker\/} of $\mathcal C$, see \cite{cornu-book}. 
Notice that the edges of $\Upsilon({\mathcal
C})$ are precisely the minimal vertex covers of $\mathcal C$. If $G$ is a
graph, then $R_s(I_c(G))$ is generated as a $K$-algebra
by elements of degree in $t$ at most two
\cite[Theorem~5.1]{cover-algebras}. One of the main result of
Section~\ref{blowupcovers1} is a graph theoretical description of the
minimal generators of 
$R_s(I_c(G))$ (see Theorem~\ref{symbo-description}). As an 
application we recover an explicit description \cite{facets}, 
in terms of
closed halfspaces, 
of the edge cone of a graph (Corollary~\ref{feasible2}). 

The symbolic
Rees algebra of the ideal $I_c({\mathcal C})$ can be interpreted in 
terms of ``$k$-vertex covers'' \cite{cover-algebras} as we now
explain.  Let $a=(a_1,\ldots,a_n)\neq 0$ be a 
vector in $\mathbb{N}^n$ and let $b\in\mathbb{N}$. We say that $a$ 
is a $b$-{\it vertex cover} of $I$ (or $\mathcal C$) 
if $\langle v_i,a\rangle\geq b$ for
$i=1,\ldots,q$. Often we will call a $b$-vertex cover simply a $b$-{\it
cover\/}. This notion plays a role in combinatorial 
optimization \cite[Chapter~77, 
p.~1378]{Schr2} and algebraic combinatorics
\cite{cover-algebras,hht-unimodular}.

The {\it algebra of covers\/} of $I$ (or $\mathcal C$), denoted by
$R_c(I)$, is the 
$K$-subalgebra of $K[t]$ generated by all monomials $x^at^b$ such 
that $a$ is a $b$-cover of $I$. We say that a $b$-cover $a$ of $I$ is
{\it reducible\/} if there exists an $i$-cover $c$ and a $j$-cover $d$ of $I$
such that $a=c+d$ and $b=i+j$. If $a$ is not reducible, we call $a$
{\it irreducible\/}. The irreducible $0$ and $1$ 
covers of $\mathcal C$ are the unit vector $e_1,\ldots,e_n$ and 
the incidence vectors $u_1,\ldots,u_s$ of the minimal vertex covers
of $\mathcal C$, respectively. The minimal generators of $R_c(I)$ as a
$K$-algebra correspond to the  
irreducible covers of $I$. Notice the
following dual descriptions:
\begin{eqnarray*}
I^{(b)}&=&(\{x^a\vert\, \langle a,u_i\rangle\geq b\mbox{ for
}i=1,\ldots,s\}), \\
J^{(b)}&=&(\{x^a\vert\, \langle a,v_i\rangle\geq b\mbox{ for
}i=1,\ldots,q\}),
\end{eqnarray*}
where $J=I_c({\mathcal C})$. Hence $R_c(I)=R_s(J)$ and $R_c(J)=R_s(I)$. 

In general each $\ell_i$ occurring in
Eq.~(\ref{okayama-car1}) determines a minimal generator of
$R_s(I_c(\mathcal{C}))$.
Indeed if we write $\ell_i=(a_i,-d_i)$, where $a_i\in\mathbb{N}^n$, 
$d_i\in\mathbb{N}$, then $a_i$
is an irreducible $d_i$-cover
of $I$ (Lemma~\ref{aug25-07}). 
Let $F_{n+1}$ be the facet of $\mathbb{R}_+(I)$ determined by the 
hyperplane $H_{e_{n+1}}$. Thus we have a map $\psi$:
$$
\begin{array}{rcl}
{\{\mbox{Facets of }\mathbb{R}_+(I({\mathcal
C}))\}\setminus\{F_{n+1}\}}&\stackrel{\psi}{\longrightarrow}&{R_s(I_c({\mathcal
C}))}\\ 
H_{\ell_k}\cap\mathbb{R}_+(I)
&\stackrel{\psi}{\longrightarrow}&{x^{a_k}t^{d_k}},
\mbox{ where }\ell_k=(a_k,-d_k)\\ 
H_{e_i}\cap\mathbb{R}_+(I)&\stackrel{\psi}{\longrightarrow}&{x_i}
\end{array}
$$
\noindent whose image provides a good approximation for the minimal set 
of generators of $R_s(I_c({\mathcal C}))$ as a $K$-algebra. Likewise the facets of 
$\mathbb{R}_+(I_c({\mathcal C}))$ give an approximation for the 
minimal set of generators of $R_s(I({\mathcal C}))$. In
Example~\ref{summer06} we show a connected graph $G$ for which the
image of 
the map $\psi$ does not generates $R_s(I_c(G))$. 
For balanced 
clutters, i.e., for clutters without odd cycles, the image of the map
$\psi$ generates $R_s(I_c({\mathcal C}))$. This follows from 
\cite[Propositions~4.10 and 4.11]{reesclu}. In particular the image of
the map $\psi$ generates $R_s(I_c({\mathcal C}))$ when $\mathcal{C}$ is a
bipartite graph. 
It would be interesting to characterize when the irreducible 
representation of the Rees cone determine the irreducible covers. 

The {\it Simis cone} of $I$ is the
rational polyhedral cone: 
$$
{\rm Cn}(I)=H_{e_1}^+\cap\cdots\cap H_{e_{n+1}}^+\cap
H_{(u_1,-1)}^+\cap\cdots\cap H_{(u_s,-1)}^+,
$$
Simis cones were
introduced in \cite{normali} to study symbolic Rees algebras of
square-free monomial ideals. If $\mathcal H$ is an integral Hilbert basis
of ${\rm Cn}(I)$, then $R_s(I(\mathcal{C}))$ equals
$K[\mathbb{N}{\mathcal H}]$, the semigroup ring of $\mathbb{N}{\mathcal H}$
(see \cite[Theorem~3.5]{normali}). This result is interesting because it 
allows us to compute the 
minimal generators of $R_s(I({\mathcal C}))$ using Hilbert bases. The
program {\it Normaliz\/} \cite{normaliz2} is suitable for computing 
Hilbert bases. There is a 
description of $\mathcal H$ valid for perfect graphs \cite{perfect}.
Perfect graphs are defined in Section~\ref{blowupcovers}.

If $G$ is a perfect graph, the irreducible
$b$-covers of $\Upsilon(G)$ 
correspond to cliques of $G$ \cite{perfect} (cf. 
Corollary~\ref{clique-description}). In 
this case, setting $\mathcal{C}=\Upsilon(G)$, 
it turns out that the image of $\psi$ generates
$R_s(I_c(\Upsilon(G)))$. Notice that $I_c(\Upsilon(G))$ is equal to
$I(G)$.  

In Section~\ref{blowupcovers} we
introduce and study the concept of an irreducible graph. A $b$-cover
$a=(a_1,\ldots,a_n)$ is called {\it binary\/} if $a_i\in\{0,1\}$ for
all $i$. We 
present a 
graph theoretical description of the irreducible binary $b$-vertex
covers of the 
blocker of $G$ (see Theorem~\ref{irr-graph-char}). It is shown that
they are in one to one  
correspondence  with the irreducible induced 
subgraphs of $G$. As a byproduct we obtain a method, using Hilbert
bases, to obtain all irreducible induced subgraphs of $G$ (see 
Corollary~\ref{method-irred-sub}). In particular we obtain all 
induced odd cycles and all induced complements of odd cycles. These
cycles are called the {\it odd holes\/} and {\it odd antiholes} of the
graph. It was shown recently \cite{fhv} that $\mathfrak{p}$ 
is an associated prime of 
$I_c(G)^2$ if and only if $\mathfrak{p}$ is generated by the 
vertices of an edge of $G$ or $\mathfrak{p}$ is generated by the
vertices of an odd hole of $G$. The proof of this remarkable 
result makes use of Theorem~\ref{symbo-description}. Odd holes and antiholes play a 
major role in graph theory. In \cite{seymour} it is shown that a graph $G$ is 
perfect if and only if $G$ has no odd holes or antiholes of length at 
least five. We give a procedure to build irreducible
graphs (Proposition~\ref{building-lemma}) and a method 
to construct irreducible $b$-vertex covers of the blocker of $G$ with
high degree relative to the number of vertices of $G$ (see
Corollaries~\ref{dec7-07} and \ref{recursive-cone}).  

Along the paper we introduce most of the 
notions that are relevant for our purposes. For unexplained
terminology we refer to \cite{diestel, Mats,Vas}. 

\section{Blowup algebras of ideals of vertex
covers}\label{blowupcovers1}

Let $G$ be a simple graph with vertex set $X=\{x_1,\ldots,x_n\}$. In
what follows we shall always assume that $G$ has no isolated vertices.
Here we will give a graph theoretical description of the irreducible
$b$-covers of $G$, i.e., we will describe the symbolic Rees 
algebra of $I_c(G)$.   

Let $S$ be a set of vertices of $G$. The {\it neighbor set\/} of $S$,
denoted by  
$N_G(S)$, is the 
set of vertices of $G$ that are adjacent with at least one 
vertex of $S$. The set $S$ is called 
{\it independent\/} if no two vertices of $S$ are adjacent. The empty
set is regarded as an independent set whose incidence vector is the
zero vector. Notice
the following duality: $S$ is a maximal 
independent set of $G$ (with respect to inclusion) if and only if
$X\setminus S$ is a minimal vertex cover of $G$.

\begin{lemma}\label{winter05} If $a=(a_i)\in\mathbb{N}^n$ is an
irreducible $k$-cover 
of $G$, then $0\leq k\leq 2$ and $0\leq a_i\leq 2$ for $i=1,\ldots,n$.
\end{lemma}

\begin{proof}Recall that $a$ is a $k$-cover of $G$ if and only if 
$a_i+a_j\geq k$ for each edge $\{x_i,x_j\}$ of
$G$. If $k=0$ or $k=1$, then by the irreducibility of $a$ 
it is seen that either $a=e_i$ for some
$i$ or 
$a=e_{i_1}+\cdots+e_{i_r}$ for some minimal vertex cover
$\{x_{i_1},\ldots,x_{i_r}\}$ of $G$. Thus we may assume that $k\geq 2$. 

Case (I): $a_i\geq 1$ for all $i$. Clearly $\mathbf{1}=(1,\ldots,1)$
is a $2$-cover. If $a-\mathbf{1}\neq 0$, then $a-\mathbf{1}$ is a
$k-2$ cover and $a=\mathbf{1}+(a-\mathbf{1})$, a contradiction 
to $a$ being an irreducible $k$-cover. Hence
$a=\mathbf{1}$. Pick any edge $\{x_i,x_j\}$ of $G$. 
Since $a$ is a $k$-cover, we get $2=a_i+a_j\geq k$ and $k$ must be
equal to $2$. 

Case (II): $a_i=0$ for some $i$. We may assume $a_i=0$ for $1\leq
i\leq r$ and $a_i\geq 1$ for $i>r$. Notice that the set
$S=\{x_1,\ldots,x_r\}$ is independent because if $\{x_i,x_j\}$ is an
edge and $1\leq i<j\leq r$, then $0=a_i+a_j\geq k$, a contradiction.
Consider the  
neighbor set $N_G(S)$ of $S$. We may assume that
$N_G(S)=\{x_{r+1},\ldots,x_s\}$. Observe that $a_i\geq k\geq 2$ for
$i=r+1,\ldots,s$, because $a$ is a $k$-cover. Write 
\begin{eqnarray*}
\lefteqn{a=(0,\ldots,0,a_{r+1}-2,\ldots,a_s-2,a_{s+1}-1,\ldots,a_n-1)+}\\
&\ \ \ \ \ \ \ \ \ \ \ \ \ \ \ \ \ \ \ \ \ \ \ \ \ \ \ \ \ \ \ \ \ \
\ \ \ \ \ \
&(\underbrace{0,\ldots,0}_r,\underbrace{2,\ldots,2}_{s-r},
\underbrace{1,\ldots,1}_{n-s})=c+d.
\end{eqnarray*}
Clearly $d$ is a $2$-cover. If $c\neq 0$, using that $a_i\geq k\geq 2$ for
$r+1\leq i\leq s$ and $a_i\geq 1$ for $i>s$ it is not hard to see that
$c$ is a $(k-2)$-cover. This gives a contradiction, because $a=c+d$.
Hence $c=0$.  
Therefore $a_i=2$ for $r<i\leq s$, $a_i=1$ for $i>s$,
and $k=2$.  \end{proof}

The next result complements the fact that the symbolic Rees algebra of
$I_c(G)$ is generated by monomials of degree in $t$ at most two 
\cite[Theorem~5.1]{cover-algebras}.

\begin{corollary} $R_s(I_c(G))$ is generated as a
$K$-algebra 
by monomials of degree in $t$ at most two and total degree at most $2n$.
\end{corollary}

\begin{proof}Let $x^at^k$ be a minimal generator of $R_s(I_c(G))$ as a
$K$-algebra. Then $a=(a_1,\ldots,a_n)$ is an irreducible $k$-cover of $G$. By
Lemma~\ref{winter05} we obtain that $0\leq k\leq 2$ and $0\leq a_i\leq
2$ for all $i$. If $k=0$ or $k=1$, we get that the degree of $x^at^k$
is at most $n$. Indeed when $k=0$ or $k=1$, one has $a=e_i$ or
$a=\sum_{x_i\in C}e_i$ for some minimal vertex cover $C$ of $G$,
respectively. If $k=2$, by the proof of Lemma~\ref{winter05} either
$a=\mathbf{1}$ 
or $a_i=0$ for some $i$. Thus $\deg(x^a)\leq 2(n-1)$. \end{proof}

Let $I=I(G)$ be the edge ideal of $G$. For use below consider the 
vectors $\ell_1,\ldots,\ell_r$ that occur in the irreducible 
representation of $\mathbb{R}_+(I)$ given in
Eq.~(\ref{okayama-car1}).  

\begin{corollary} If
$\ell_i=(\ell_{i1},\ldots,\ell_{in},-\ell_{i(n+1)})$, then $0\leq
\ell_{ij}\leq 2$ for $j=1,\ldots,n$ and $1\leq\ell_{i(n+1)}\leq 2$. 
 \end{corollary}

\begin{proof}It suffices to observe that $(\ell_{i1},\ldots,\ell_{in})$ is 
an irreducible $\ell_{i(n+1)}$-cover of $G$ and to apply 
Lemma~\ref{winter05}. \end{proof}

\begin{lemma}\label{jan15-05} $a=(1,\ldots,1)$ is an irreducible
$2$-cover of $G$ if  
and only if $G$ is non bipartite.  
\end{lemma}

\begin{proof}$\Rightarrow$) We proceed by contradiction assuming that $G$
is a bipartite graph. Then $G$ has a bipartition $(V_1,V_2)$. Set 
$a'=\sum_{x_i\in V_1}e_i$ and $a''=\sum_{x_i\in V_2}e_i$. Since
$V_1$ and $V_2$ are minimal vertex covers of $G$, we can decompose 
$a$ as $a=a'+a''$, where $a'$ and $a''$ are $1$-covers, which is
impossible. 

$\Leftarrow$) Notice that $a$ cannot be the sum of a $0$-cover and a
$2$-cover. Indeed if $a=a'+a''$, where $a'$ is a $0$-cover and $a''$
is a $2$-cover, then $a''$ has an entry $a_i$ equal to zero. Pick 
an edge $\{x_i,x_j\}$ incident with $x_i$, then $\langle
a'',e_i+e_j\rangle\leq 1$, a contradiction. 
Thus we may assume that $a=c+d$, where $c,d$ are $1$-covers. Let
$C_r$ be an odd 
cycle of $G$ of length $r$. Notice that any vertex cover of $C_r$
must contain 
a pair of adjacent vertices because $r$ is odd. Clearly a vertex cover
of $G$ is also a vertex cover of the subgraph $C_r$. Hence the vertex covers
of $G$ corresponding to $c$ and $d$ must contain a pair of adjacent
vertices, a contradiction because $c$ and $d$ are complementary
vectors and the complement of a vertex cover is an 
independent set. \end{proof}

\begin{definition}\label{bala-clu}\rm  
Let $A$ be the incidence matrix of a clutter $\mathcal{C}$. A clutter
$\mathcal{C}$ has a {\it cycle\/} of length $r$ if there is a square
sub-matrix of $A$ of 
order $r\geq 3$ with exactly two $1$'s in 
each row and column.  A clutter without odd cycles is 
called {\it balanced}.
\end{definition}

\begin{proposition}{\rm(\cite[Proposition~4.11]{reesclu})}\label{balanced-ntf} 
If $\mathcal{C}$ is a balanced clutter, then 
$$
R_s(I_c(\mathcal{C}))=R[I_c(\mathcal{C})t].
$$
\end{proposition} 

This result was first shown for bipartite graphs in 
\cite[Corollary~2.6]{alexdual} and later generalized to balanced
clutters \cite{reesclu} using an algebro combinatorial description 
of clutters with the max-flow min-cut property \cite{clutters}. 

Let $S$ be a set of vertices of a graph $G$, the {\it induced
subgraph\/} on $S$, denoted by $\langle S\rangle$, is the maximal subgraph of $G$ with
vertex set $S$. The next result has been used in \cite{fhv} to show
that any associated prime of $I_c(G)^2$ is generated by the 
vertices of an edge of $G$ or it is generated by the
vertices of an odd hole of $G$. 

We come to the main result of this section.

\begin{theorem}\label{symbo-description} 
Let $0\neq a=(a_i)\in\mathbb{N}^n$ and let $\Upsilon(G)$ be the family
of minimal vertex covers of a graph $G$. 
\begin{enumerate}
\item[\rm(i)\ ] If $G$ is bipartite, then $a$ is an irreducible
$b$-cover of $G$ 
if and only if\/ $b=0$ and $a=e_i$ for some $1\leq i\leq n$ or $b=1$ and
$a=\sum_{x_i\in
C}e_i$ for some $C\in\Upsilon(G)$.
\item[\rm (ii)] If $G$ is non-bipartite, then $a$ is an irreducible
$b$-cover if and only if $a$ has one of the following forms{\rm :} 
\begin{enumerate}
\item[\rm (a)] $($$0$-covers$)$ $b=0$ and $a=e_i$ for some $1\leq i\leq n$,

\item[\rm (b)] $($$1$-covers$)$ $b=1$ and $a=\sum_{x_i\in C}e_i$ for
some $C\in\Upsilon(G)$,

\item[\rm (c)] $($$2$-covers$)$ $b=2$ and $a=(1,\ldots,1)$,

\item[\rm (d)] $($$2$-covers$)$ $b=2$ and up to permutation of vertices 
$$a=(\underbrace{0,\ldots,0}_{|A|},\underbrace{2,\ldots,2}_{|N_G(A)|},
{1,\ldots,1})$$
for some independent set of vertices $A\neq\emptyset$ of $G$ such
that
\begin{enumerate}
\item[$(\mathrm{d}_1)$] 
$N_G(A)$ is not a vertex cover of $G$ and $V\neq A\cup N_G(A)$, 
\item[$(\mathrm{d}_2)$] the induced subgraph $\langle
V\setminus(A\cup N_G(A))\rangle$ has no isolated vertices 
and is not bipartite.
\end{enumerate}
\end{enumerate}
\end{enumerate}
\end{theorem}

\begin{proof}(i) $\Rightarrow$) Since $G$ is bipartite,  by
Proposition~\ref{balanced-ntf}, we have 
the equality $R_s(J)=R[Jt]$,
where $J=I_c(G)$ is the ideal of vertex covers of $G$. Thus the 
minimal  set of generator of $R_s(J)$ as a $K$-algebra is 
the set 
$$
\{x_1,\ldots,x_n,x^{u_1}t,\ldots,x^{u_s}t\},$$ 
where $u_1,\ldots,u_s$ are the incidence vectors of the minimal
vertex covers of $G$. By hypothesis $a$ is an irreducible $b$-cover of
$G$, i.e., $x^at^b$ is a minimal generator of 
$R_s(I_c(\mathcal{C}))$. Therefore either $a=e_i$ for some $i$ and
$b=0$ or 
$a=u_i$ for some $i$ and $b=1$. The
converse follows readily and is valid for any graph or clutter.  

(ii) $\Rightarrow)$ By Lemma~\ref{winter05} we have $0\leq b\leq 2$ 
and $0\leq a_i\leq 2$ for all $i$. If $b=0$ or $b=1$, then clearly $a$
has the form indicated in (a) or (b) respectively. 

Assume $b=2$. If $a_i\geq 1$ for all $i$, the $a_i=1$ for all $i$,
otherwise if $a_i=2$ for some $i$, then $a-e_i$ is a $2$-cover and 
$a=e_i+(a-e_i)$, a contradiction. Hence $a=\mathbf{1}$. Thus we may
assume that $a$ has the form 
$$
a=(0,\ldots,0,2,\ldots,2,1,\ldots,1).
$$
We set $A=\{x_i\vert\, a_i=0\}\neq\emptyset$, $B=\{x_i\vert\,
a_i=2\}$, 
and $C=V\setminus(A\cup B)$. Observe that $A$ is an independent set
because $a$ is a $2$-cover and $B=N_G(A)$ because $a$ 
is irreducible. Hence it is seen that conditions $(\mathrm{d}_1)$ 
and $(\mathrm{d}_2)$ are satisfied. Using Lemma~\ref{jan15-05}, the
proof of the converse is direct. \end{proof}

\begin{lemma}\label{aug25-07} Let $\mathcal{C}$ be a clutter and let
$I=I(\mathcal{C})$ be its edge ideal. If $\ell_k=(a_k,-d_k)$ is any
of the vectors that occur 
in Eq.~{\rm(\ref{okayama-car1})}, where $a_k\in\mathbb{N}^n$, 
$d_k\in\mathbb{N}$, then $a_k$ is an irreducible $d_k$-cover of
$\mathcal{C}$. 
\end{lemma}

\begin{proof} We proceed by contradiction assume there is a $d_k'$-cover
$a_k'$ and a $d_k''$-cover $a_k''$ such that $a_k=a_k'+a_k''$ and
$d_k=d_k'+d_k''$. Set $F'=H_{(a_k',-d_k')}\cap \mathbb{R}_+(I)$ and 
$F''=H_{(a_k'',-d_k'')}\cap \mathbb{R}_+(I)$. Clearly $F',F''$ are
proper faces of $\mathbb{R}_+(I)$ and  $F=\mathbb{R}_+(I)\cap
H_{\ell_k}=F'\cap F''$. Recall that any proper face of
$\mathbb{R}_+(I)$ is the intersection of those facets that contain it
(see \cite[Theorem~3.2.1(vii)]{webster}). Applying this fact to 
$F'$ and $F''$ it is seen that $F'\subset F$ or $F''\subset F$, i.e.,
$F=F'$ or $F=F''$. We may assume $F=F'$. Hence 
$H_{(a_k',-d_k')}=H_{\ell_k}$. Taking orthogonal complements 
we get that $(a_k',-d_k')=\lambda(a_k,-d_k)$ for some 
$\lambda\in\mathbb{Q}_+$, because the orthogonal complement of 
$H_{\ell_k}$ is one dimensional and it is generated by $\ell_k$. Since the non-zero 
entries of $\ell_k$ are relatively prime, we may assume that
$\lambda\in\mathbb{N}$. Thus $d_k'=\lambda d_k\geq d_k\geq d_k'$ and
$\lambda$ must be equal to $1$. Hence $a_k=a_k'$ and $a_k''$ must be
zero, a contradiction. \end{proof}

\begin{example}\label{summer06}\rm Consider the following graph $G$:
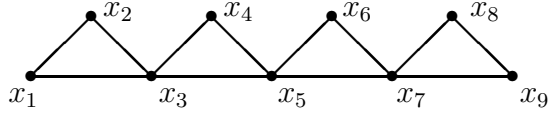
\begin{center}
$
\setlength{\unitlength}{.02cm}
\thicklines  
\begin{picture}(100,50) 
\put(-150,0){\circle*{6.4}}
\put(-150,0){\line(1,0){320}}
\put(-30,40){\circle*{6.4}}
\put(50,40){\circle*{6.4}}
\put(130,40){\circle*{6.4}}
\put(-110,40){\circle*{6.4}}
\put(-102,40){$x_2$}
\put(-22,40){$x_4$}
\put(57,40){$x_6$}
\put(142,40){$x_8$}
\put(-70,0){\circle*{6.4}}
\put(-70,0){\line(1,1){40}}
\put(10,0){\line(1,1){40}}
\put(90,0){\line(1,1){40}}
\put(10,0){\circle*{6.4}}
\put(90,0){\circle*{6.4}}
\put(170,0){\circle*{6.4}}
\put(-165,-17){$x_1$}
\put(-65,-17){$x_3$}
\put(14,-17){$x_5$}
\put(93,-17){$x_7$}
\put(175,-17){$x_9$}
\put(-110,40){\line(1,-1){40}}
\put(-110,40){\line(-1,-1){40}}
\put(-30,40){\line(1,-1){40}}
\put(50,40){\line(1,-1){40}}
\put(130,40){\line(1,-1){40}}
\end{picture}
\vspace{0.5cm}
$
\end{center}
Using {\it Normaliz\/} \cite{normaliz2} it is seen that the vector
$a=(1,1,2,0,2,1,1,1,1)$ is an irreducible $2$-cover of $G$ such that
the supporting hyperplane $H_{(a,-2)}$ does not define a facet of the
Rees cone of $I(G)$. Thus, in general, the image of $\psi$ described
in the introduction does not 
determine $R_s(I_c(G))$. We may use Lemma~\ref{jan15-05} to
construct non-connected graphs with this property. 
\end{example}

\paragraph{Edge cones of graphs} Let $G$ be a connected simple graph and let
${\mathcal A}=\{v_1,\ldots,v_q\}$ be the set of all vectors $e_i+e_j$
such that 
$\{x_i,x_j\}$ is an 
edge of $G$. The {\em edge cone\/}  of $G$, 
denoted by ${\mathbb R}_+{\mathcal A}$, is defined as the 
cone  generated by ${\mathcal A}$. Below we give an explicit
combinatorial  
description of the edge cone. 

Let $A$ be an {\it independent set\/} of vertices of $G$. The
supporting 
hyperplane of the edge cone of $G$ defined by  
\[
\sum_{x_i\in N_G(A)}x_i-\sum_{x_i\in A}x_i=0
\]
will be denoted by $H_A$. 

Edge cones and their representations by closed halfspaces are a
useful tool to study $a$-invariants of edge 
subrings \cite{shiftcon,join}. The following result is a prototype 
of these representations. As an application we give a 
direct proof of the next result using Rees cones.

\begin{corollary}{\rm(\cite[Corollary~2.8]{facets})}\label{feasible2}
A  vector $a=(a_1,\ldots,a_n)\in {\mathbb R}^n$  
is in ${\mathbb R}_+{\mathcal A}$ if and only if $a$ satisfies the
following system of linear inequalities 
\[
\begin{array}{rcll}
a_i&\geq & 0, &i=1,\ldots,n;
\\ & & &\vspace{-3mm} \\ 
\sum_{x_i\in N_G(A)}a_i-\sum_{x_i\in A}a_i& \geq & 0,& \mbox{for
all independent sets }A\subset V(G).
\end{array}
\] 
\end{corollary}

\begin{proof}Set ${\mathcal B}=\{(v_1,1),\ldots,(v_q,1)\}$ and $I=I(G)$. Notice
the equality
\begin{equation}\label{ottawa}
\mathbb{R}_+(I)\cap\mathbb{R}{\mathcal B}=\mathbb{R}_+{\mathcal B},
\end{equation}
where $\mathbb{R}{\mathcal B}$ is $\mathbb{R}$-vector space spanned by
$\mathcal B$. 
Consider the irreducible representation of $\mathbb{R}_+(I)$ given 
in Eq.~(\ref{okayama-car1}) and write $\ell_i=(a_i,-d_i)$, where
$0\neq a_i\in\mathbb{N}^n$, $0\neq d_i\in\mathbb{N}$. Next we show the
equality: 
\begin{equation}\label{sept5-1-07}
{\mathbb R}_+{\mathcal A}={\mathbb R}{\mathcal A}\cap\mathbb{R}_+^n\cap
H_{(2a_1/d_1-\mathbf{1})}^+ \cap\cdots\cap
H_{(2a_r/d_r-\mathbf{1})}^+,
\end{equation}
where $\mathbf{1}=(1,\ldots,1)$. Take $\alpha\in {\mathbb R}_+{\mathcal
A}$. Clearly $\alpha\in\mathbb{R}\mathcal{A}\cap \mathbb{R}^n_+$. We can write
$$ 
\alpha=\lambda_1v_1+\cdots+\lambda_qv_q\ \Rightarrow\
|\alpha|=2(\lambda_1+\cdots+\lambda_q)=2b.
$$
Thus $(\alpha,b)=\lambda_1(v_1,1)+\cdots+\lambda_q(v_q,1)$, i.e., 
$(\alpha,b)\in\mathbb{R}_+\mathcal{B}$. Hence from 
Eq.~(\ref{ottawa}) we get $(\alpha,b)\in\mathbb{R}_+(I)$ and
$$
\langle(\alpha,b),(a_i,-d_i) \rangle\geq 0\ \Rightarrow\ 
\langle\alpha,a_i\rangle\geq bd_i=(|\alpha|/2)d_i=|\alpha|(d_i/2). 
$$
Writing $\alpha=(\alpha_1,\ldots,\alpha_n)$ and
$a_i=(a_{i1},\ldots,a_{in})$, the last inequality gives:
$$
\alpha_1a_{i_1}+\cdots+\alpha_na_{in}\geq
(\alpha_1+\cdots+\alpha_n)(d_i/2)\ \Rightarrow\ 
\langle\alpha,a_i-({d_i}/{2})\mathbf{1}\rangle\geq 0. 
$$
Then $\langle\alpha,2a_i/d_i-\mathbf{1}\rangle\geq 0$ and 
$\alpha\in H_{(2a_i/d_i-\mathbf{1})}^+$ for all $i$, as required.
This proves that 
$\mathbb{R}_+{\mathcal A}$ is contained in the right hand side of
Eq.~(\ref{sept5-1-07}). The other inclusion follows similarly. Now by 
Lemma~\ref{aug25-07} we obtain that $a_i$ is an irreducible
$d_i$-cover of $G$. Therefore, using the explicit description 
of the irreducible $b$-covers of $G$ given in
Theorem~\ref{symbo-description}, we get the equality
$$
{\mathbb R}_+{\mathcal A}=\left(\bigcap_{A\in \mathcal F} H_A^+\right)
\bigcap\left(\bigcap_{i=1}^n H_{e_i}^+\right),
$$
where $\mathcal F$ is the collection of all the independent sets of
vertices  of $G$. From this equality the assertion follows at
once. \end{proof}

The edge cone of $G$ encodes information about both the Hilbert
function of the  
edge subring $K[G]$ (see \cite{shiftcon}) and the graph $G$ itself. 
As a simple illustration, we recover the following version 
of the marriage theorem for bipartite graphs, see \cite{Boll}. Recall
that a pairing by an independent set of edges of all the vertices of
a graph $G$ is called a  
{\it perfect matching\/} or a $1$-{\it factor}. 

\begin{corollary} If $G$ is a bipartite connected graph, then  $G$ has
a perfect 
matching if and only 
if $|A|\leq |N_G(A)|$ for every independent
set of vertices $A$ of $G$.
\end{corollary}

\begin{proof}Notice that the graph $G$ has a perfect matching if and only if the 
vector $\beta=(1,1,\ldots,1)$ is in $\mathbb{N}{\mathcal A}$. By 
\cite[Lemma~2.9]{shiftcon} we have the equality
$\mathbb{Z}^n\cap\mathbb{R}_+{\mathcal 
A}=\mathbb{N}{\mathcal A}$. Hence $\beta$ is in $\mathbb{N}{\mathcal A}$ if 
and only if $\beta\in\mathbb{R}_+{\mathcal A}$. Thus the result 
follows from Corollary~\ref{feasible2}. \end{proof}

\section{Symbolic Rees algebras of  edge ideals}\label{blowupcovers}

Let $G$ be a graph with vertex set $X=\{x_1,\ldots,x_n\}$ and let 
$I=I(G)$ be its edge ideal. As before we denote the clutter of
minimal vertex covers of $G$ by $\Upsilon(G)$. The clutter
$\Upsilon(G)$ is called the {\it blocker\/} of $G$. Recall that the
symbolic Rees algebra of $I(G)$ is given by
\begin{equation}\label{sept5-07}
R_s(I(G))=K[x^{a}t^b\vert\, a\mbox{ is an irreducible }b\mbox{-cover
of }  \Upsilon(G)],
\end{equation}
where the set $\{x^{a}t^b\vert\, a\mbox{ is an irreducible }b\mbox{-cover
of }  \Upsilon(G)\}$ is the minimal set of generators of $R_s(I(G))$
as a $K$-algebra. The main purpose of
this section is to study the symbolic Rees 
algebra of $I(G)$ via graph theory. We are interested in finding
combinatorial 
representations for the minimal set of generators of this algebra. 

\begin{lemma}\label{july4-06} Let $0\neq
a=(a_1,\ldots,a_m,0,\ldots,0)\in\mathbb{N}^n$ and  
let $a'=(a_1,\ldots,a_m)$. If $0\neq b\in\mathbb{N}$, then $a$ is an
irreducible $b$-cover of 
$\Upsilon(G)$  
if and only if $a'$ is an irreducible $b$-cover of 
$\Upsilon(\langle S\rangle)$, where $S=\{x_1,\ldots,x_m\}$.
\end{lemma}

\begin{proof}It suffices to prove that $a$ is a $b$-cover of the blocker of
$G$ if and only if $a'$ is a $b$-cover of the blocker of $\langle
S\rangle$. 

$\Rightarrow)$ The induced subgraph $\langle
S\rangle$ is not a discrete graph. 
Take a minimal vertex cover $C'$ of $\langle S\rangle$.
Set $C=C'\cup(V(G)\setminus S)$. Since $C$ is a vertex cover of $G$ 
such that $C\setminus\{x_i\}$ is not a vertex cover of $G$ for every
$x_i\in C'$,
there is a minimal vertex cover $C_\ell$ of $G$ such that $C'\subset
C_\ell\subset C$ and $C'=C_\ell\cap S$. Notice that
$$
\sum_{x_i\in C'}a_i=\sum_{x_i\in C_\ell\cap S}a_i=\langle
a,u_\ell\rangle\geq b,
$$
where $u_\ell$ is the incidence vector of $C_\ell$. 
Hence $\sum_{x_i\in C'}a_i\geq b$, as required.

$\Leftarrow$) Take a minimal vertex cover $C_\ell$ of $G$. Then
$C_\ell\cap S$ contains a minimal vertex cover $C_\ell'$ of $\langle
S\rangle$. Let $u_\ell$ (resp. $u_\ell'$) be the incidence vector of 
$C_\ell$ (resp. $C_\ell'$). Notice that
$$
\langle a,u_\ell\rangle=\sum_{x_i\in C_\ell\cap S}a_i
\geq \sum_{x_i\in C_\ell'}a_i=\langle
a',u_\ell'\rangle\geq b.
$$
Hence $\langle a,u_\ell\rangle\geq b$, as required. \end{proof}

We denote a complete subgraph of $G$ with $r$ vertices by ${\mathcal
K}_r$. If $v$ is a vertex of $G$, we denote its neighbor set by
$N_G(v)$.  

\begin{lemma}\label{july4-06-1} Let $G$ be a graph and let
$a=(a_1,\ldots,a_n)$ be an 
irreducible 
$b$-cover of
$\Upsilon(G)$ such that $a_i\geq 1$ for all $i$. If 
$\langle N_G(x_n)\rangle={\mathcal K}_r$, then $a_i=1$ for all $i$, 
$b=r$, $n=r+1$, and $G={\mathcal K}_{n}$.
\end{lemma}

\begin{proof}We may assume that $N_G(x_n)=\{x_1,\ldots,x_r\}$. 
We set 
$$
c=e_1+\cdots+e_r+e_n; \ \ \
d=(a_1-1,\ldots,a_r-1,a_{r+1},\ldots,a_{n-1},a_n-1). 
$$
Notice that $\langle x_1,\ldots,x_{r},x_n\rangle={\mathcal K}_{r+1}$.
Thus $c$ is an $r$-cover 
of $\Upsilon(G)$ because any minimal vertex cover of $G$ must 
intersect all edges of ${\mathcal K}_{r+1}$. By the irreducibility of $a$,
there exists a minimal vertex cover $C_\ell$ of $G$ such that 
$\sum_{x_i\in C_\ell}a_i=b$. Clearly we have $b\geq g\geq r$, where $g$ is
the height of $I(G)$. Let $C_k$ be an arbitrary minimal vertex cover
of $G$. Since $C_k$ contains exactly $r$ vertices of ${\mathcal K}_{r+1}$,
from the inequality $\sum_{x_i\in C_k}a_i\geq b$ we get 
$\sum_{x_i\in C_k}d_i\geq b-r$, where $d_1,\ldots,d_n$ are the entries
of $d$. Therefore $d=0$; otherwise if $d\neq 0$, then $d$ is a
$(b-r)$-cover of $\Upsilon(G)$ and $a=c+d$, a contradiction to the
irreducibility of $a$.   
It follows that $g=r$, $n=r+1$, $a_i=1$ 
for $1\leq i\leq r$, $a_n=1$, and $G={\mathcal K}_{n}$. \end{proof}

\medskip

\noindent {\it Notation} We regard ${\mathcal K}_0$ as
the empty set with zero elements. A sum over an empty set is
defined to be $0$. 

\begin{proposition}\label{myperfect-char-mod} Let $G$ be a graph and let 
$J=I_c(G)$ be its ideal of vertex covers. Then the set 
$$
F=\{(a_i)\in\mathbb{R}^{n+1}\vert\,
\textstyle\sum_{x_i\in\mathcal{K}_r}a_i=(r-1)a_{n+1}\}\cap\mathbb{R}_+(J)
$$
is a facet of the Rees cone $\mathbb{R}_+(J)$. 
\end{proposition}

\begin{proof}If $\mathcal{K}_r=\emptyset$, then
$r=0$ and $F=H_{e_{n+1}}\cap\mathbb{R}_+(J)$, which is clearly a
facet because  
$e_1,\ldots,e_n\in F$. If $r=1$, then $F=H_{e_i}\cap\mathbb{R}_+(J)$ 
for some $1\leq i\leq n$, which is a facet because $e_j\in F$ for
$j\notin\{i,n+1\}$ and there is at least 
one minimal vertex cover of $G$ not containing $x_i$. We may assume that
$X'=\{x_1,\ldots,x_r\}$ is the vertex set of
${\mathcal K}_r$ and $r\geq 2$. For each
$1\leq i\leq r$ there is a minimal vertex cover $C_i$ of $G$ not 
containing $x_i$. Notice that $C_i$ contains $X'\setminus\{x_i\}$.
Let $u_i$ be the incidence vector 
of $C_i$. Since the rank of $u_1,\ldots,u_r$ is $r$, it follows that
the set 
$$
\{(u_1,1),\ldots,(u_r,1),e_{r+1},\ldots,e_n\}
$$
is a linearly independent set contained in $F$, i.e., $\dim(F)=n$.
Hence $F$ is a facet of $\mathbb{R}_+(J)$ because the hyperplane that
defines $F$ is a supporting hyperplane. \end{proof}

\begin{proposition}\label{cliques-are-irr} Let $G$ be a graph and let 
$0\neq a=(a_i)\in\mathbb{N}^n$. If 
\begin{enumerate}
\item[\rm (a)] $a_i\in\{0,1\}$ for all $i$, and 

\item[\rm (b)] $\langle\{x_i\vert\, a_i>0\}\rangle={\mathcal K}_{r+1}$,
\end{enumerate}
then $a$ is an irreducible $r$-cover of $\Upsilon(G)$.
\end{proposition}

\begin{proof}By Proposition~\ref{myperfect-char-mod}, 
the closed halfspace $H_{(a,-r)}^+$ occurs in the irreducible representation 
of the Rees cone $\mathbb{R}_+(J)$, where $J=I_c(G)$. Hence $a$ is
an irreducible $r$-cover by Lemma~\ref{aug25-07}. \end{proof}

A {\it clique\/} of a graph $G$ is a set of vertices that induces
a complete subgraph. We will also call a complete subgraph of $G$ 
a clique. Symbolic Rees algebras are related to perfect graphs as is 
seen below.  Let us recall the notion of perfect graph. 
A {\it colouring\/} of the vertices of $G$ is an assignment
of colours to the vertices of $G$ in such a way that adjacent vertices
have distinct colours. The {\it chromatic number\/} of $G$ 
is the minimal number of colours in a colouring of $G$.
A graph is {\it perfect\/} if for every induced subgraph $H$, the
chromatic 
number of $H$ equals the size of the largest complete subgraph 
of $H$. We refer 
to \cite {cornu-book,diestel,Schr2} and the references there 
for the theory of perfect graphs.

\medskip

\noindent {\it Notation\/} The {\it support\/} of 
$x^a=x_1^{a_1}\cdots x_n^{a_n}$ is ${\rm supp}(x^a)= \{x_i\, |\,
a_i>0\}$. 
\begin{corollary}[\rm\cite{perfect}]\label{clique-description}
If $G$ is a graph, then 
$$
K[x^at^r\vert\, 
x^a\mbox{ square-free };\, 
\langle{\rm supp}(x^a)\rangle={\mathcal K}_{r+1};\, 0\leq
r<n]\subset R_s(I(G))
$$
with equality if and only if $G$ is a perfect graph. 
\end{corollary}

\begin{proof}The inclusion follows from Proposition~\ref{cliques-are-irr}. 
If $G$ is a perfect graph, then by \cite[Corollary~3.3]{perfect} the
equality holds. Conversely if the equality holds, then by
Lemma~\ref{aug25-07} and Proposition~\ref{myperfect-char-mod} we have 
\begin{equation}\label{mundial06}
\mathbb{R}_+(I_c(G))=\left\{(a_i)\in\mathbb{R}^{n+1}\vert\,
\textstyle\sum_{x_i\in\mathcal{K}_r}a_i
\geq (r-1)a_{n+1};\ \forall\,  {\mathcal K}_r\subset
G\right\}.
\end{equation}
Hence a direct application of \cite[Proposition~2.2]{perfect} 
gives that $G$ is a perfect graph.   \end{proof}

The {\it vertex
covering number\/} of $G$, denoted by $ \alpha_0 (G)$, is the number of
vertices in a minimum vertex cover of $G$ (the cardinality of any smallest
vertex cover in $G$). Notice that $\alpha_0(G)$ equals the height of
$I(G)$. If $H$ is a discrete graph, i.e., all the vertices of $H$ are
isolated, we set $I(H)=0$ and $\alpha_0(H)=0$. 

\begin{lemma}\label{rats} Let $G$ be a graph. If
$a=e_1+\cdots+e_r$ is an irreducible $b$-cover of 
$\Upsilon(G)$, then $b=\alpha_0(H)$, 
where $H=\langle x_1,\ldots,x_r\rangle$.
\end{lemma}

\begin{proof}The case $b=0$ is clear. Assume $b\geq 1$. Let $C_1,\ldots,C_s$
be the minimal vertex covers of $G$ and 
let $u_1,\ldots,u_s$ be their incidence vectors. Notice that $\langle
a,u_i\rangle=b$ for some $i$. Indeed if $\langle
a,u_i\rangle> b$ for all $i$, then $a-e_1$ is a $b$-cover of
$\Upsilon(G)$ and $a=(a-e_1)+e_1$, a contradiction. Hence
$$
b=\langle a,u_i\rangle=|\{x_1,\ldots,x_r\}\cap C_i|\geq\alpha_0(H).
$$
This proves that $b\geq\alpha_0(H)$. Notice that $H$ is not a discrete
graph. Then we can pick a minimal
vertex cover $A$ of $H$ such that $|A|=\alpha_0(H)$. The set 
$$
C=A\cup(V(G)\setminus\{x_1,\ldots,x_r\})
$$ 
is a vertex cover of $G$. Hence there is a minimal vertex 
cover $C_\ell$ of $G$ such that $A\subset C_\ell\subset C$. Observe
that $C_\ell\cap \{x_1,\ldots,x_r\}=A$. Thus we
get $\langle a,u_\ell\rangle=|A|\geq b$, i.e., $\alpha_0(H)\geq b$.
Altogether we have $b=\alpha_0(H)$. \end{proof}

This result has been recently extended to clutters using the notion
of parallelization \cite{cm-mfmc}.  
Let $\mathcal{C}$ be a clutter on the vertex set
$X=\{x_1,\ldots,x_n\}$ and let $x_i\in X$. Then {\it duplicating\/}
$x_i$ means extending $X$ by a new vertex $x_i'$ and replacing
$E(\mathcal{C})$ by
$$    
E(\mathcal{C})\cup\{(e\setminus\{x_i\})\cup\{x_i'\}\vert\, x_i\in e\in
E(\mathcal{C})\}.
$$
The {\it deletion\/} of $x_i$, denoted by
$\mathcal{C}\setminus\{x_i\}$, is the clutter formed from
$\mathcal{C}$ by deleting the vertex $x_i$ and all edges containing
$x_i$. A clutter obtained from $\mathcal{C}$ by a sequence of
deletions and 
duplications of vertices is called a {\it parallelization\/}. If 
$w=(w_i)$ is a vector in $\mathbb{N}^n$, we denote by $\mathcal{C}^w$
the clutter obtained from
$\mathcal{C}$ by deleting any vertex $x_i$ with $w_i=0$ and
duplicating $w_i-1$ times any vertex $x_i$ if $w_i\geq 1$. The map
$w\mapsto \mathcal{C}^w$ gives a one to one correspondence between
$\mathbb{N}^n$ and the parallelizations of $\mathcal{C}$. 

\begin{example}
Let $G$ be the graph whose only edge is 
$\{x_1,x_2\}$ and let $w=(3,3)$. Then
$G^w=\mathcal{K}_{3,3}$ is the complete bipartite graph 
with bipartition $V_1=\{x_1,x_1^2,x_1^3\}$ and
$V_2=\{x_2,x_2^2,x_2^3\}$. Notice that $x_i^k$ is a vertex, i.e., $k$
is an index not an exponent.
\end{example}

\begin{proposition}[\rm \cite{cm-mfmc}] Let $\mathcal{C}$ be a clutter and let 
$\Upsilon(\mathcal{C})$ be the blocker of 
$\mathcal{C}$. If $w=(w_i)$ is an irreducible $b$-cover of
$\Upsilon(\mathcal{C})$, then
$$
b=\left.\min\left\{\sum_{x_i\in
C}w_i\right\vert\,
C\in\Upsilon(\mathcal{C})\right\}
= \alpha_0(\mathcal{C}^w).
$$
\end{proposition}

The next result gives a nice graph theoretical description for 
the irreducible binary $b$-vertex
covers of the 
blocker of $G$.

\begin{theorem}\label{irr-graph-char} Let $G$ be a graph and let $a=(1,\ldots,1)$. 
Then $a$ is a reducible $\alpha_0(G)$-cover of $\Upsilon(G)$ if 
and only if there are $H_1$ and $H_2$ induced subgraphs of $G$ such that
\begin{enumerate}
\item[\rm (i)\ ] $V(G)$ is the disjoint union of $V(H_1)$ and $V(H_2)$, and 
\item[\rm (ii)] $\alpha_0(G)=\alpha_0(H_1)+\alpha_0(H_2)$.   
\end{enumerate}
\end{theorem}

\begin{proof}$\Rightarrow$) We may assume that $a_1=e_1+\cdots+e_r$, 
$a_2=a-a_1$, $a_i$ is a $b_i$-cover of $\Upsilon(G)$, $b_i\geq 1$ for
$i=1,2$, and 
$\alpha_0(G)=b_1+b_2$. Consider the subgraphs $H_1=\langle
x_1,\ldots,x_r\rangle$  and $H_2=\langle x_{r+1},\ldots,x_n\rangle$.
Let $A$ be a minimal vertex cover of $H_1$ with $\alpha_0(H_1)$
vertices. Since 
$$C=A\cup(V(G)\setminus\{x_1,\ldots,x_r\})$$
is a
vertex cover $G$, there is a minimal vertex cover $C_k$ of $G$ such
that $A\subset C_k\subset C$. Hence 
$$
|A|=|C_k\cap\{x_1,\ldots,x_r\}|=\langle
a_1,u_k\rangle\geq b_1,
$$  
and $\alpha_0(H_1)\geq b_1$. Using a similar argument we get that 
$\alpha_0(H_2)\geq b_2$. If $C_\ell$ is a minimal vertex
cover of $G$ with $\alpha_0(G)$ vertices, then $C_\ell\cap V(H_i)$ is a
vertex cover of $H_i$. Therefore
$$
b_1+b_2=\alpha_0(G)=|C_\ell|=\sum_{i=1}^2|C_\ell\cap V(H_i)|\geq 
\sum_{i=1}^2\alpha_0(H_i)\geq b_1+b_2,
$$
and consequently $\alpha_0(G)=\alpha_0(H_1)+\alpha_0(H_2)$. 

$\Leftarrow$) We may assume that $V(H_1)=\{x_1,\ldots,x_r\}$ and 
$V(H_2)=V(G)\setminus V(H_1)$. Set $a_1=e_1+\cdots+e_r$ and 
$a_2=a-a_1$. For any minimal vertex cover $C_k$ of $G$, we have that 
$C_k\cap V(H_i)$ is a vertex cover of $H_i$. Hence 
$$
\langle
a_1,u_k\rangle=|C_k\cap\{x_1,\ldots,x_r\}|\geq\alpha_0(H_1),
$$
where $u_k$ is the incidence vector of $C_k$. Consequently $a_1$ is
an $\alpha_0(H_1)$-cover  of $\Upsilon(G)$. Similarly
we obtain that $a_2$ is an $\alpha_0(H_2)$-cover of $\Upsilon(G)$. Therefore
$a$ is a reducible $\alpha_0(G)$-cover of $\Upsilon(G)$. \end{proof}
 
\begin{definition}\rm A graph satisfying conditions (i) and (ii) is
called a {\it reducible} graph. If $G$ is not reducible, it is called
{\it irreducible}.
\end{definition}

These notions appear in \cite{erdos-gallai}. As far as we know there is no
structure theorem for irreducible graphs. Examples of
irreducible graphs include complete graphs, odd cycles, and
complements of odd cycles. Below we give a method, using Hilbert
bases, to obtain all irreducible induced subgraphs of $G$. 

\medskip

By \cite[Lemma~5.4]{korte} there exists a finite set 
${\mathcal H}\subset \mathbb{N}^{n+1}$ such 
that 
\begin{enumerate}
\item[\rm (a)] ${\rm Cn}(I(G))=\mathbb{R}_+{\mathcal H}$, 
and 
\item[\rm (b)] $\mathbb{Z}^{n+1}\cap
\mathbb{R}_+{\mathcal H}=\mathbb{N}{\mathcal H}$, 
\end{enumerate}
where $\mathbb{N}{\mathcal H}$ is the additive subsemigroup of 
$\mathbb{N}^{n+1}$ generated 
by ${\mathcal H}$. 
\begin{definition}\label{hilb-basis-def}
\rm The set ${\mathcal H}$ is called a 
{\it Hilbert basis\/} of ${\rm Cn}(I(G))$. 
\end{definition}

If we 
require $\mathcal{H}$ 
to be minimal (with respect inclusion), then 
$\mathcal{H}$ is unique \cite{Schr1}.

\begin{corollary}\label{method-irred-sub} Let $G$ be a graph and let
$\alpha=(a_1,\ldots,a_n,b)$ be a vector in $\{0,1\}^n\times\mathbb{N}$. Then 
$\alpha$ is an element of the minimal integral Hilbert basis of ${\rm
Cn}(I(G))$ if and only if the induced subgraph 
$H=\langle\{x_i\vert\, a_i=1\}\rangle$ is irreducible with
$b=\alpha_0(H)$.
\end{corollary}

\begin{proof}The map $(a_1,\ldots,a_n,b)\mapsto x_1^{a_1}\cdots
x_n^{a_n}t^b$ establishes a  
one to one correspondence between the minimal integral Hilbert basis of ${\rm
Cn}(I(G))$ and the minimal generators of $R_s(I(G))$ as 
a $K$-algebra. Thus the result follows from 
Lemma~\ref{july4-06} and Theorem~\ref{irr-graph-char}. \end{proof}

The next result shows that irreducible
graphs occur naturally in the theory of perfect graphs. 

\begin{proposition} 
A graph $G$ is perfect if and only if the only irreducible induced
subgraphs of $G$ are the complete subgraphs.
\end{proposition} 

\begin{proof}$\Rightarrow$) Let $H$ be an irreducible induced subgraph of
$G$. We may assume that $V(H)=\{x_1,\ldots,x_r\}$. Set
$a'=(1,\ldots,1)\in\mathbb{N}^r$ and
$a=(a',0\ldots,0)\in\mathbb{N}^n$. By Theorem~\ref{irr-graph-char},
$a'$ is an irreducible
$\alpha_0(H)$ cover of $\Upsilon(H)$. Then by Lemma~\ref{july4-06}, 
$a'$ is an irreducible $\alpha_0(H)$ cover of $\Upsilon(G)$. Since $x_1\cdots
x_rt^{\alpha_0(H)}$ is a minimal generator of $R_s(I(G))$, using
Corollary \ref{clique-description}  we obtain that $\alpha_0(H)=r-1$
and that $H$ is a complete subgraph 
of $G$ on $r$ vertices. 

$\Leftarrow$) In \cite{seymour} it is shown that $G$ is a perfect graph 
if and only if no induced subgraph of $G$ is an odd cycle of length at 
least five or the complement of one. Since odd cycles 
and their complements are irreducible subgraphs. It follows that $G$
is perfect. \end{proof}

\begin{definition}\rm A graph $G$ is called {\it vertex critical\/} if 
$\alpha_0{(G\setminus\{x_i\})}<\alpha_0{(G)}$ for all $x_i\in V(G)$.          
\end{definition}

\begin{remark}\label{nov29-07}\rm If $x_i$ is any vertex of a graph $G$ and 
$\alpha_0{(G\setminus\{x_i\})}<\alpha_0{(G)}$, then
$\alpha_0(G\setminus\{x_i\})=\alpha_0(G)-1$
\end{remark}

\begin{lemma} If the graph $G$ is irreducible, then it is connected 
and vertex critical  
\end{lemma}

\begin{proof}Let $G_1,\ldots,G_r$ be the connected components of $G$. Since 
$\alpha_0(G)$ is equal to $\sum_{i}\alpha_0(G_i)$, we get $r=1$. Thus
$G$ is connected. To complete the proof it suffices to prove that 
$\alpha_0(G\setminus\{x_i\})<\alpha_0(G)$ for all $i$ (see
Remark~\ref{nov29-07}). If
$\alpha_0(G\setminus\{x_i\})=\alpha_0(G)$, then
$G=H_1\cup H_2$, where $H_1=G\setminus\{x_i\}$ and $V(H_2)=\{x_i\}$, 
a contradiction. \end{proof}

\begin{definition}\rm The {\it cone\/} $C(G)$, over a graph $G$, 
is obtained by adding a new vertex $v$ to $G$ and joining every vertex
of $G$ to $v$.
\end{definition}

The next result can be used to build irreducible graphs. In particular
it follows that cones over irreducible graphs are irreducible.

\begin{proposition}\label{building-lemma} Let $G$ be a graph with $n$
vertices and let $H$ be a graph obtained from $G$ by adding a new
vertex $v$ and some new edges joining $v$ with $V(G)$. 
If $a=(1,\ldots,1)\in\mathbb{N}^n$ is an irreducible
$\alpha_0(G)$-cover  of $\Upsilon(G)$ such that
$\alpha_0(H)=\alpha_0(G)+1$, then $a'=(a,1)$
is an irreducible $\alpha_0(H)$-cover of $\Upsilon(H)$. 
\end{proposition}

\begin{proof}Clearly $a'$ is an $\alpha_0(H)$-cover of $\Upsilon(H)$. 
Assume that $a'=a_1'+a_2'$, where $0\neq a_i'$ is a $b_i'$-cover of
$\Upsilon(H)$ and $b_1'+b_2'=\alpha_0(H)$. We may assume that
$a_1'=(1,\ldots,1,0,\ldots,0)$ and $a_2'=(0,\ldots,0,1,\ldots,1)$. Let
$a_i$ be the vector in $\mathbb{N}^n$ obtained from $a_i'$ by 
removing its last entry. Set $v=x_{n+1}$. Take a minimal vertex cover
$C_k$ of $G$ 
and consider $C_k'=C_k\cup\{x_{n+1}\}$. Let $u_k'$ (resp. $u_k$) be 
the incidence vector of $C_k'$ (resp. $C_k$). Then
$$
\langle a_1,u_k\rangle=\langle a_1',u_k'\rangle\geq b_1'\ \mbox{and}\ 
\langle a_2,u_k\rangle+1=\langle a_2',u_k'\rangle\geq b_2',
$$
consequently $a_1$ is a $b_1'$-cover of $\Upsilon(G)$. If
$b_2'=0$, then $a_1$ is an $\alpha_0(H)$-cover of $\Upsilon(G)$, 
a contradiction; because if $u$ is the incidence vector of a minimal
vertex cover of $G$ with $\alpha_0(G)$ elements, then we would obtain
$\alpha_0(G)\geq\langle u,a_1\rangle\geq\alpha_0(H)$, which is
impossible. Thus $b_2'\geq 1$, and $a_2$ is a $(b_2'-1)$-cover of
$\Upsilon(G)$ if $a_2\neq 0$.
Hence $a_2=0$, because $a=a_1+a_2$ and $a$ is 
irreducible. This means that $a_2'=e_{n+1}$ is a $b_2'$-cover of
$\Upsilon(H)$, a contradiction. Therefore $a'$ is an irreducible
$\alpha_0(H)$-cover of $\Upsilon(H)$, as required. \end{proof}

\begin{definition}\rm A graph $G$ is called {\it edge critical\/} if 
$\alpha_0{(G\setminus e)}<\alpha_0{(G)}$ for all $e\in E(G)$.           
\end{definition}

\begin{proposition} If $G$ is a connected edge critical graph, 
then $G$ is irreducible.
\end{proposition}

\begin{proof}Assume that $G$ is reducible. Then there are induced subgraphs 
$H_1$, $H_2$ of $G$ such that $V(H_1)$, $V(H_2)$ is a partition of 
$V(G)$ and $\alpha_0(G)=\alpha_0(H_1)+\alpha_0(H_2)$. Since $G$ is
connected there is an edge $e=\{x_i,x_j\}$ with 
$x_i$ a vertex of $H_1$ and $x_j$ a vertex of $H_2$. Pick a minimal vertex cover $C$ 
of $G\setminus e$ with $\alpha_0(G)-1$ vertices. As $E(H_i)$ is a 
subset of $E(G\setminus e)=E(G)\setminus\{e\}$ for $i=1,2$, we get that 
$C$ covers all edges of $H_i$ for $i=1,2$. Hence $C$ must have at 
least $\alpha_0(G)$ elements, a contradiction. \end{proof}

\begin{corollary} The following hold for any connected
graph{\rm :}
$$
\begin{array}{cccccc}
\mbox{edge critical}&\Longrightarrow& \mbox{irreducible}&
&\Longrightarrow
&\mbox{vertex critical.}
\end{array}
$$
\end{corollary}

\paragraph{Finding generators of symbolic Rees algebras using cones}

The {\em cone\/} $C(G)$, over the graph $G$, is obtained by adding a new vertex $t$ 
to $G$ and joining every vertex of $G$ to $t$. 

\begin{example}\rm A pentagon and its cone:

\vspace{1.3cm}

$
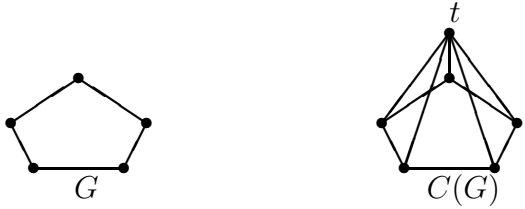
\begin{array}{cccc}
\ \ \ \ \ \ \ \ \ \ \ \ \ \ \ \ \ \ \ \ \ \ \ &
\setlength{\unitlength}{.03cm}
\thicklines
\begin{picture}(40,30)
\put(-30,0){\circle*{4.2}}
\put(30,0){\circle*{4.2}}
\put(0,20){\circle*{4.2}}
\put(-20,-20){\circle*{4.2}}
\put(20,-20){\circle*{4.2}}
\put(-2,-33){$G$}
\put(-30,0){\line(3,2){30}}
\put(30,0){\line(-3,2){30}}
\put(-30,0){\line(1,-2){10}}
\put(-20,-20){\line(1,0){40}}
\put(20,-20){\line(1,2){10}}
\end{picture}
&
\ \ \ \ \ \ \ \ \ \ \ \ \ \ \ \ \ \ \ \ \ \ 
&
\setlength{\unitlength}{.03cm}
\thicklines
\begin{picture}(40,30)
\put(-30,0){\circle*{4.2}}
\put(30,0){\circle*{4.2}}
\put(0,20){\circle*{4.2}}
\put(-20,-20){\circle*{4.2}}
\put(20,-20){\circle*{4.2}}
\put(0,40){\circle*{4.2}}
\put(0,45){$t$}
\put(-10,-33){$C(G)$}
\put(30,0){\line(-3,2){30}}
\put(-30,0){\line(3,2){30}}
\put(-30,0){\line(1,-2){10}}
\put(-20,-20){\line(1,0){40}}
\put(20,-20){\line(1,2){10}}
\put(0,40){\line(0,-1){20}}
\put(0,40){\line(-3,-4){30}}
\put(0,40){\line(3,-4){30}}
\put(0,40){\line(-1,-3){20}}
\put(0,40){\line(1,-3){20}}
\end{picture}
\end{array}
$
\end{example}

\vspace{1.1cm}

In \cite{bahiano} Bahiano showed that if $H=C(G)$ is the graph
obtained by taking a 
cone over a pentagon $G$ with vertices $x_1,\ldots,x_5$, then 
$$
R_s(I(H))=R[I(H)t][x_1\cdots x_5t^3,x_1\cdots x_6t^4,x_1\cdots
x_5x_6^2t^5]. 
$$
This simple example shows that taking a cone over an irreducible graph
tends to increase the degree in $t$ of the generators of the 
symbolic Rees algebra. Other examples using this ``cone process''
have been 
shown in \cite[Example~5.5]{cover-algebras}. 

Let $G$ be a graph with vertex set $V(G)=\{x_1,\ldots,x_n\}$. 
The aim here is to give 
a general procedure---based on the irreducible representation of the
Rees cone of $I_c(G)$---to construct generators of $R_s(I(H))$ of 
high degree in $t$, where $H$ is a graph constructed from $G$ by
recursively taking cones over graphs already constructed. 

By the finite basis theorem
\cite[Theorem~4.11]{webster} there is a 
unique irreducible representation
\begin{equation}\label{okayama-dual-alex} 
{\mathbb R}_+(I_c(G))=H_{e_1}^+\cap H_{e_2}^+\cap\cdots\cap
H_{e_{n+1}}^+\cap H_{\alpha_1}^+\cap H_{\alpha_2}^+\cap\cdots\cap
H_{\alpha_p}^+
\end{equation}
such that each $\alpha_k$ is in $\mathbb{Z}^{n+1}$, the non-zero
entries of 
each $\alpha_k$ are relatively prime, and none of the closed 
halfspaces $H_{e_1}^+,\ldots,
H_{e_{n+1}}^+,H_{\alpha_1}^+,\ldots,H_{\alpha_p}^+$ can be
omitted from the intersection. For use below we assume that 
$\alpha$ is any of the vectors $\alpha_1,\ldots,\alpha_p$
that occur in the irreducible representation. Thus we can write 
$\alpha=(a_1,\ldots,a_n,-b)$ for some $a_i\in\mathbb{N}$ and for 
some $b\in\mathbb{N}$.

\begin{lemma}\label{nov30-07} Let $H$ be the cone over $G$. If
$$
\beta=(a_1,\ldots,a_n,\textstyle
(\sum_{i=1}^na_i)-b,\textstyle
-\sum_{i=1}^na_i)=(\beta_1,\ldots,\beta_{n+1},-\beta_{n+2})
$$ and 
$a_i\geq 1$ for all $i$, then $F=H_\beta\cap \mathbb{R}_+(I_c(H))$
is a facet of $\mathbb{R}_+(I_c(H))$.
\end{lemma}

\begin{proof}First we prove that $\mathbb{R}_+(I_c(H))\subset
H_\beta^{+}$, 
i.e., $H_\beta$ is a supporting hyperplane of the Rees cone. 
By Lemma~\ref{aug25-07}, $(a_1,\ldots,a_n)$ is an irreducible
$b$-cover of $\Upsilon(G)$. Hence there is $C\in\Upsilon(G)$ such 
that $\sum_{x_i\in C}a_i=b$. Therefore $\beta_{n+1}$ 
is greater or equal than $1$. This proves that
$e_1,\ldots,e_{n+1}$ are in $H_\beta^+$. Let $C$ be any 
minimal vertex cover of $H$ and let $u=\sum_{x_i\in C}e_i$ 
be its characteristic vector. Case (i): If $x_{n+1}\notin C$, 
then $C=\{x_1,\ldots,x_n\}$ and 
$$
\sum_{x_i\in C}\beta_i=\sum_{i=1}^na_i=\beta_{n+2},
$$
that is, $(u,1)\in H_\beta^+$. Case (ii): If $x_{n+1}\in C$, then
$C_1=C\setminus\{x_{n+1}\}$ is a minimal vertex cover of $G$. Hence 
$$
\sum_{x_i\in C}\beta_i=\sum_{x_i\in C_1}\beta_i+\beta_{n+1}\geq 
b+\beta_{n+1}=\beta_{n+2},
$$
that is, $(u,1)\in H_\beta^+$. Therefore $\mathbb{R}_+(I_c(H))\subset
H_\beta^{+}$. To prove that $F$ is a facet we must show it has
dimension $n+1$ because the dimension of $\mathbb{R}_+(I_c(H))$ is
$n+2$. We denote the characteristic vector of a minimal vertex cover
$C_k$ of $G$ by $u_k$. 
By hypothesis there are minimal vertex covers $C_1,\ldots,C_n$
of $G$ such that the vectors $(u_1,1),\ldots,(u_n,1)$ are linearly
independent and  
\begin{equation}\label{dec1-07}
\langle (a,-b),(u_k,1)\rangle=0\ \Longleftrightarrow  \langle
a,u_k\rangle=b,
\end{equation}
for $k=1,\ldots,n$. Therefore
\begin{eqnarray*}
\langle(\beta_1,\ldots,\beta_{n+1}),(u_k,1)\rangle=\beta_{n+2}\ \mbox{
and }\ \ \ \ \ \ \ \ \ \ \ \ \ \ \ \ \ & &\\ 
\langle(\beta_1,\ldots,\beta_{n+1}),(1,\ldots,1,0)\rangle=\beta_{n+2},
&
&
\end{eqnarray*}
i.e., the set $\mathcal{B}=\{(u_1,1),\ldots,(u_n,1),(1,\ldots,1,0\}$ 
is contained in $H_\beta$. Since
$$
C_1\cup\{x_{n+1}\},\ldots,C_n\cup\{x_{n+1}\},\{x_1,\ldots,x_n\}$$
are minimal vertex covers of $H$, the set $\mathcal{B}$ is also
contained in $\mathbb{R}_+(I_c(H))$ and consequently in $F$. 
Thus its suffices to prove that
$\mathcal{B}$ is linearly independent. If $(1,\ldots,1,0)$ is a linear
combination of $(u_1,1),\ldots,(u_n,1)$, then we can write
$$
(1,\ldots,1)=\lambda_1u_1+\cdots+\lambda_nu_n
$$
for some scalars $\lambda_1,\ldots,\lambda_n$ such that
$\sum_{i=1}^n\lambda_i=0$. Hence from Eq.~(\ref{dec1-07}) we get
$$
|a|=\langle (1,\ldots,1),a\rangle=\lambda_1\langle u_1,a\rangle
+\cdots+\lambda_n\langle u_n,a\rangle=(\lambda_1+\cdots+\lambda_n)b=0,
$$
a contradiction. \end{proof}
\begin{corollary}\label{dec7-07} If $a_i\geq 1$ for all $i$, 
then $x_1^{\beta_1}\cdots x_{n+1}^{\beta_{n+1}}t^{\beta_{n+2}}$ is 
a minimal generator of $R_s(I(H))$.
\end{corollary}

\begin{proof}By Lemma~\ref{nov30-07}, $F=H_\beta\cap \mathbb{R}_+(I_c(H))$
is a facet of $\mathbb{R}_+(I_c(H))$. Therefore using 
Lemma~\ref{aug25-07}, the vector $(\beta_1,\ldots,\beta_{n+1})$ is an 
irreducible $\beta_{n+2}$-cover of $\Upsilon(H)$, i.e., 
$x_1^{\beta_1}\cdots x_{n+1}^{\beta_{n+1}}t^{\beta_{n+2}}$ is 
a minimal generator of $R_s(I(H))$.  \end{proof}

\begin{corollary}\label{recursive-cone} Let $G_0=G$ and let $G_r$ be
the cone over  
$G_{r-1}$ for $r\geq 1$. If $\alpha=(1,\ldots,1,-g)$, then 
$$
(\underbrace{1,\ldots,1}_{n},\underbrace{n-g,\ldots,n-g}_{r})
$$
is an irreducible $n+(r-1)(n-g)$ cover of $G_r$. In particular 
$R_s(I(G_r))$ has a generator of degree in $t$ equal to 
$n+(r-1)(n-g)$.
\end{corollary}

As a very particular example of our construction consider: 

\begin{example}\rm Let $G=C_s$ be an odd cycle of length $s=2k+1$. Note
that $\alpha_0(C_s)=(s+1)/2=k+1$. Then 
by Corollary~\ref{recursive-cone} 
$$
x_1\cdots x_sx_{s+1}^k\cdots x_{s+r}^kt^{rk+k+1}
$$
is a minimal generator of $R_s(I(G_r))$. This illustrates that the degree in
$t$ of the minimal generators of $R_s(I(G_r))$ is much larger than 
the number of vertices of the graph $G_r$ \cite{cover-algebras}. 
\end{example}

\end{document}